\documentclass[a4pitaper,11pt]{amsart}
\usepackage{amsmath,amsthm,amssymb}
\usepackage[mathscr]{eucal}
 \usepackage{cite}
\usepackage{upgreek}
\usepackage[bookmarks,bookmarksnumbered]{hyperref}
\usepackage{enumerate}

\numberwithin{equation}{section}

\newcommand{\car}{\curvearrowright}

\theoremstyle{plain}
\newtheorem{main}{Theorem}
\newtheorem{mcor}[main]{Corollary}

\newtheorem{theorem}{Theorem}[section]

\newtheorem{lemma}[theorem]{Lemma}
\newtheorem{proposition}[theorem]{Proposition}

\theoremstyle{definition}

\newtheorem{definition}[theorem]{Definition}

\newtheorem{notation}[theorem]{Notation}
\newtheorem{remark}[theorem]{Remark}

\begin{document}


\title[Product rigidity in von Neumann and C$^*$-algebras via s-malleable deformations]
{Product rigidity in von Neumann and C$^*$-algebras via s-malleable deformations}

\author{Daniel Drimbe}
\address{Department of Mathematics, KU Leuven, Celestijnenlaan 200b, B-3001 Leuven, Belgium}
\email{daniel.drimbe@kuleuven.be}
\thanks {The author holds the postdoctoral fellowship fundamental research 12T5221N of the Research Foundation - Flanders.}

\begin{abstract} 
We provide a new large class of countable icc groups $\mathcal A$ for which the product rigidity result from \cite{CdSS15} holds: if $\Gamma_1,\dots,\Gamma_n\in\mathcal A$ and $\Lambda$ is any group such that $L(\Gamma_1\times\dots\times\Gamma_n)\cong L(\Lambda)$, then there exists a product decomposition $\Lambda=\Lambda_1\times\dots\times \Lambda_n$ such that $L(\Lambda_i)$ is stably isomorphic to $L(\Gamma_i)$, for any $1\leq i\leq n$. Class $\mathcal A$ consists of groups $\Gamma$ for which $L(\Gamma)$ admits an s-malleable deformation in the sense of Sorin Popa and it includes all non-amenable groups $\Gamma$ such that either (a) $\Gamma$ admits an unbounded 1-cocycle into its left regular representation, or (b) $\Gamma$ is an arbitrary wreath product group with amenable base. As a byproduct of these results, we obtain new examples of W$^*$-superrigid groups and new rigidity results in the C$^*$-algebra theory.
\end{abstract}

\maketitle

\section{Introduction}
Every countable group $\Gamma$ gives rise to the group von Neumann algebra $L(\Gamma)$ by considering the weak operator closure of the complex group algebra $\mathbb C[\Gamma]$ acting on the Hilbert space $\ell^2(\Gamma)$ by left convolution \cite{MvN43}. A main theme in operator algebras is the classification of group von Neumann algebras which is centered around the following question: what properties of the group $\Gamma$ are remembered by 
$L(\Gamma)$? This problem is the most interesting when $\Gamma$ is icc (i.e., all non-trivial conjugacy classes of $\Gamma$ are infinite), which corresponds to $L(\Gamma)$ being a II$_1$ factor. In the amenable case, the classification is completed by the work of Connes \cite{Co76} which asserts that any two icc amenable groups give rise to isomorphic von Neumann algebras. Therefore, besides the amenability of the group, no information can be recovered from $L(\Gamma)$ when $\Gamma$ is icc amenable.

In the non-amenable case the situation is radically different and far more complex. An outstanding progress has been achieved since the invention of Popa’s deformation/rigidity theory \cite{Po07} and there have been discovered many instances when various algebraic and analytical properties of a group $\Gamma$ can be recovered from $L(\Gamma)$, see the surveys \cite{Va10a, Io12b,Io17}. Remarkably, Ioana, Popa and Vaes found in \cite{IPV10} the first class of countable groups $\Gamma$ that are W$^*$-superrigid. Roughly speaking, this means that the group $\Gamma$ is completely remembered by its von Neumann algebra $L(\Gamma)$. Subsequently, several other classes of W$^*$-superrigid groups have been found \cite{BV12,Be14,CI17,CD-AD20}.



However, in general, one can only expect to recover certain aspects of a group $\Gamma$ from its von Neumann algebra $L(\Gamma)$. We only highlight the following developments. Ozawa showed in \cite{Oz03} that the group von Neumann algebra of a non-amenable bi-exact icc group is prime, in particular implying that $L(\Gamma)\ncong L(\Gamma_1\times\Gamma_2)$, for all infinite groups $\Gamma_1$ and $\Gamma_2$. Ozawa and Popa then proved that any tensor product of II$_1$ factors of non-amenable hyperbolic icc groups admits a unique prime decomposition into prime factors \cite{OP03}. As a corollary, their result shows that if  $L(\Gamma_1\times\dots\times\Gamma_n)\cong L(\Lambda_1\times\dots\times\Lambda_m)$ for some icc hyperbolic groups $\Gamma_i$'s and infinite groups $\Lambda_j$'s, then $m=n$ and after a permutation of indices we have $L(\Gamma_i)$ is stably isomorphic to $L(\Lambda_i)$, for all $1\leq i\leq n$. More recently, Chifan, de Santiago and Sinclair strengthened the previous corollary of \cite{OP03} by discovering the following {\it product rigidity} phenomenon: if the groups $\Gamma_i$'s are icc hyperbolic, then any group $\Lambda$ such that 
$L(\Gamma_1\times\dots\times\Gamma_n)\cong L(\Lambda)$ admits a decomposition $\Lambda=\Lambda_1\times\dots\times\Lambda_n$ satisfying $L(\Gamma_i)$ is stably isomorphic to $L(\Lambda_i)$, for any $1\leq i\leq n$ \cite{CdSS15}. While a large number of other unique prime factorization results have been obtained since \cite{OP03} (see, e.g., the introduction of \cite{DHI16}), the above product rigidity result has been extended to the class of non-amenable bi-exact groups only very recent \cite{CD-AD20}.
Some of the methods from \cite{CD-AD20}, including an augmentation technique, will play an important role in our work as well.

The goal of this paper is to provide a new class of countable groups, denoted Class $\mathcal A$, for which the above product rigidity holds, see Theorem \ref{A}. A common feature of these groups is that their von Neumann algebras admit an s-malleable deformation in the sense of Popa \cite{Po01,Po03} (see Definition \ref{def:malleable}) and a key ingredient of the proof of our first main result is the use of Popa's spectral gap principle that was developed in \cite{Po06a,Po06b}. In fact, Theorem \ref{A} follows from a more general and conceptual result, see Theorem \ref{Th:general}.

{\bf Class $\mathcal{A}$.} We say that a countable non-amenable icc group $\Gamma$ belongs to Class $\mathcal A$ if $\Gamma$ satisfies one of the following conditions:
\begin{enumerate}
    \item $\Gamma$ admits an unbounded cocycle for some mixing representation $\pi:\Gamma\to \mathcal O( H_{\mathbb R})$ such that $\pi$ is weakly contained in the left regular representation of $\Gamma$.
    
    \item $\Gamma=\Gamma_1*_\Sigma\Gamma_2$ is an amalgamated free product group satisfying $[\Gamma_1:\Sigma]\ge 2$ and $[\Gamma_2:\Sigma]\ge 3$, where $\Sigma<\Gamma$ is an amenable almost malnormal\footnote{A subgroup $H<G$ is called almost malnormal if $gHg^{-1}\cap H$ is finite for any $g\in G\setminus$H.} subgroup.

    \item $\Gamma=\Sigma\wr_{G/H} G$ is a generalized wreath product group with $\Sigma$ amenable, $G$ non-amenable and $H<G$ is an amenable almost malnormal subgroup. 
\end{enumerate}

For any $i\in\overline{1,3}$, if $\Gamma$ belongs to $\mathcal A$ satisfying condition (i), then we say that $\Gamma$ belongs to $\mathcal A_i$. Note that $\mathcal A_1$ contains all non-amenable icc groups $\Gamma$ with $\beta^{(2)}_1(\Gamma)>0$.






\begin{main}\label{A}
Let $\Gamma=\Gamma_1\times\dots\times\Gamma_n$ be a product of $n\ge 1$ countable groups that belong to $\mathcal A$ and denote $M=L(\Gamma)$. Let $\Lambda$ be any countable group such that $M^t=L(\Lambda)$ for some $t>0.$

Then there exist a product decomposition $\Lambda=\Lambda_1\times\dots\times\Lambda_n$, a unitary $u\in\mathcal U(M^t)$ and some positive numbers $t_1,\dots,t_n$ with $t_1\dots t_n=t$ such that $uL(\Lambda_i)u^*=L(\Gamma_i)^{t_i}$, for any $i\in\overline{1,n}.$
\end{main}

Remark that \cite[Theorem A]{CdSS15} shows that if the groups $\Gamma_i$'s are non-amenable free groups and $L(\Gamma_1\times\dots\times\Gamma_n)\cong L(\Lambda)$ for some group $\Lambda$, then $\Lambda$ admits a product decomposition into $n$ infinite groups. Theorem \ref{A} strengthens this fact by replacing the groups $\Gamma_i$'s with the more general class of arbitrary free product groups.


To put our result into a better perspective, we note that Popa's deformation/rigidity theory led to striking rigidity results for group von Neumann algebras of wreath product groups. To recall this results, fix a non-trivial abelian group $A.$ Popa showed in \cite{Po03,Po04} that the group von Neumann algebras $L(A\wr\Gamma_1)$ are pairwise non-isomorphic for different icc property (T) groups $\Gamma_1$. This result was strengthened by Ioana, Popa and Vaes in \cite{IPV10}  by showing that for any icc property (T) group $\Gamma_1$, the isomorphism $L(A\wr\Gamma_1)\cong L(\Lambda)$ implies that there exists a semi-direct product decomposition $\Lambda=B\rtimes\Lambda_1$ such that $\Gamma_1\cong\Lambda_1.$ Several other rigidity results have been obtained for group von Neumann algebras of (generalized) wreath product groups, including primeness, unique prime factorization and relative solidity, see \cite{Io06,Po06a,CI08,IPV10,SW11,BV12,IM19}. Theorem \ref{A} provides on the other hand a new general rigidity result for wreath product groups by showing that the von Neumann algebra of a product of wreath product groups with amenable base completely remembers the product structure.

\begin{remark}\label{newrigidity}
We also emphasize the following new rigidity phenomenon that Theorem \ref{A} (more precisely, Theorem \ref{Th:general}) leads to. First, we consider the class $\mathscr M_0$ of all non-amenable II$_1$ factors $M$ that admit an s-malleable deformation $(\tilde M, (\alpha_t)_{t\in\mathbb R})$ with the properties that $M\subset\tilde M$ is mixing, $_ML^2(\tilde M)\ominus L^2(M)_M$ is weakly contained in the coarse bimodule $_ML^2(M)\otimes L^2(M)_M$ and $\alpha_t$ does not converge uniformly on $(M)_1$ (see also Subsection \ref{class M}). In Theorem \ref{Th:general} we classify all tensor product decompositions of any group von Neumann algebra for which the tensor factors belong to $\mathscr{M}_0$. 
In contrast to the product rigidity result from \cite{CdSS15}, the II$_1$ factors from $\mathscr{M}_0$ are not necessarily group von Neumann algebras. For instance, any tracial non-amenable free product $M_1*M_2$ belongs to $\mathscr M_0$, see Remark \ref{inclusion}.
\end{remark}




%

Next, we show that Theorem \ref{A} can be used together with \cite{IPV10} to derive new examples of W$^*$-superrigid product groups. First, recall that a countable group $\Gamma$ is {\it W$^*$-superrigid} if for any group $\Lambda$ and any $*$-isomorphism $\theta:L(\Gamma)^t\to L(\Lambda)$ for some $t>0$, we have $t=1$ and there exist a group isomorphism $\delta:\Gamma\to\Lambda$, a character $\omega:\Gamma\to\mathbb T$ and a unitary $w\in L(\Lambda)$ such that $\theta(u_g)=\omega(g)w v_{\delta(g)}w^*$, for any $g\in\Gamma$. Here, we denoted by $\{u_g\}_{g\in\Gamma}$ and $\{v_\lambda\}_{\lambda\in\Lambda}$ the canonical generating unitaries of $L(\Gamma)$ and $L(\Lambda)$, respectively.



The first class of W$^*$-superrigid product groups has recently been found in \cite{CD-AD20} by considering products of W$^*$-superrigid groups from \cite{IPV10} that
are bi-exact. As a consequence of Theorem \ref{A} we can actually drop the bi-exactness assumption and therefore obtain that all products of W$^*$-superrigid groups from \cite{IPV10} are again W$^*$-superrigid. To illustrate our result, we introduce the following class of groups that was considered in \cite{IPV10}.

{\bf Class $\mathcal{IPV}$.} We say that a countable group $\Gamma$ belongs to Class $\mathcal {IPV}$ if $\Gamma=(\mathbb Z/n\mathbb Z)\wr_{I} G$ is a generalized wreath product group that satisfies:
\begin{itemize}
    \item $n\in\{2,3\}$ and $I=G/H$, where $H<G$ is an infinite amenable almost malnormal subgroup.
    \item $G$ admits an infinite normal subgroup that either has relative property (T) or its centralizer is non-amenable.
    \end{itemize}



\begin{mcor}\label{B}
If $\Gamma=\Gamma_1\times\dots\times\Gamma_n$ is a product of  W$^*$-superrigid groups that belong to $\mathcal A$ (e.g., $\Gamma_i\in \mathcal {IPV}$ for any $i$), then $\Gamma$ is W$^*$-superrigid.
\end{mcor}

The problem of proving that the W$^*$-superrigid property is closed with respect to direct products is notoriously hard and remains open. In Corollary \ref{B} (see also Theorem \ref{Th:general}) we make some progress on this problem, by showing that within the class of generalized wreath product groups with almost malnormal stabilizers, the W$^*$-superrigidity property is preserved by taking direct products.


Finally, we will discuss some applications of Theorem \ref{A} to the C$^*$-algebra theory. In contrast to the von Neumann algebra setting, the classification of reduced $C^*$-algebras is not governed by an amenable/non-amenable dichotomy in the sense that the reduced group C$^*$-algebra $C_r^*(\Gamma)$ of an amenable group $\Gamma$ does not provide the same striking lack of rigidity as its von Neumann algebra $L(\Gamma)$. In fact, any torsion free abelian group $\Gamma$ is $C^*$-superrigid \cite[Theorem 1]{Sc74}, which roughly means that $C_r^*(\Gamma)$ completely remembers the group $\Gamma$. 
By building upon the result of \cite{Sc74}, several other classes of amenable C$^*$-superrigid groups have been found, see the introduction of \cite{ER18}. In the non-amenable case, the only examples of C$^*$-superrigid groups are obtained in \cite{CI17,CD-AD20} via Popa's deformation/rigidity theory using their von Neumann algebraic superrigid behavior combined with the unique trace property \cite{BKKO14} of their reduced C$^*$-algebra. In a similar way, the von Neumann rigidity result from Theorem \ref{A} can be transferred to a product rigidity in $C^*$-algebras.

\begin{mcor}\label{C}
Let $\Gamma=\Gamma_1\times\dots\times\Gamma_n$ be a product of $n\ge 1$ countable non-amenable icc groups such that $\beta^{(2)}_1(\Gamma_i)>0$ or $\Gamma_i\in \mathcal A_2$, for any $i\in\overline{1,n}$. Let $\Lambda$ be any countable group satisfying $C^*_r(\Gamma)=C^*_r(\Lambda)$.

Then $\Lambda=\Lambda_1\times\dots\times\Lambda_n$ admits a product decomposition into infinite groups such that $L(\Lambda_i)$ is stably isomorphic to $L(\Gamma_i)$, for any $i\in\overline{1,n}$.
\end{mcor}


We note that the above product rigidity has already been obtained for icc hyperbolic groups \cite{CdSS15} and extended to infinite direct sums of icc hyperbolic property (T) groups \cite{CU18}. Finally, it would be interesting to show that in certain cases the conclusion of Corollary \ref{C} can be strengthened to derive that the C$^*$-algebras of the factors $C^*_r(\Gamma_i)$ and $C^*_r(\Lambda_i)$ are isomorphic for any $i\in\overline{1,n}$.


{\bf Outline of the proof of Theorem \ref{A}.} We outline briefly and informally the proof of our main result, which is Theorem \ref{A}. Let $\Gamma=\Gamma_1\times\dots\times\Gamma_n$ be a product of $n\ge 1$ countable groups that belong to Class $\mathcal A$; to better illustrate  the proof, we assume that $\beta^{(2)}_1(\Gamma_i)>0$, for any $i\in\overline{1,n}$. Denote $M_i=L(\Gamma_i)$, for any $i\in \overline{1,n}$ and let $M=L(\Gamma).$

Our goal is to show that for any group von Neumann algebra decomposition $M^t=L(\Lambda)$, where $t>0$ and $\Lambda$ is a countable group, the underlying group is a product group $\Lambda=\Lambda_1\times\dots \times\Lambda_n$. To simplify notation, we assume $t=1$. In order to attain this goal, we use a number of techniques from Popa's deformation/rigidity theory. Following \cite{PV09}, we define the $*$-homomorphism $\Delta:M\to M\bar\otimes M$ by letting $\Delta(v_\lambda)=v_\lambda\otimes v_\lambda$, for any $\lambda\in\Lambda$.

Throughout the proof we repeatedly use  Popa's spectral gap rigidity principle together with the s-malleable deformation constructed by Sinclair \cite{Si10} for all $L(\Gamma_i)$'s. 
Their first combined application shows that
\begin{equation}\label{intro1}
    \forall i\in\overline{1,n}, \; \exists t_i\in\overline{1,n} \text{\; \; such that \; \; } \Delta(M_{\widehat i})\prec M\bar\otimes M_{\widehat {t_i}}.
\end{equation}
Here, $P\prec Q$ denotes that a corner of $P$ embeds into $Q$ inside the ambient algebra, in the sense of Popa \cite{Po03} and $M_{\widehat i}=\bar\otimes_{j\neq i}M_j$.
Next, by combining Ioana's ultrapower technique \cite{Io11} with some recent techniques in the framework of von Neumann algebras without property Gamma \cite{BMO19,IM19}, we are able to upgrade \eqref{intro1} to
\begin{equation}\label{intro2}
    \Delta(M_I)\prec M\bar\otimes M_{I}, \text{\;\; for any }I\subset \overline{1,n},
\end{equation}
We also derive that for any subset $I\subset\overline{1,n}$ there exists a subgroup $\Sigma_I<\Lambda$ with non-amenable centralizer $C_\Lambda(\Sigma)$ such that
\begin{equation}\label{intro3}
L(\Sigma_I)\prec M_{I} \text{\; \; and \; \;} M_I\prec L(\Sigma_I).
\end{equation}
Here, we denoted by $M_I=\bar\otimes_{i\in I}M_i.$ This is achieved in Section 4. Next, by using some general recent results on s-malleable deformations \cite{dSHHS20}
and
an augmentation technique developed in \cite{CD-AD20} we are able to use \eqref{intro2} and \eqref{intro3} to show that there exists a non-zero projection $e\in L(\Sigma_{\widehat n})'\cap M$ such that 
\begin{equation}\label{intro4}
    L(\Sigma_{\widehat n})e\prec^s M_{\hat n} \text{\; \; and \; \; } M_{\hat n}\prec L(\Sigma_{\widehat n})e. 
\end{equation}
Here, we denoted by $P\prec^s Q$ to indicate that $Pp'\prec Q$ for any non-zero projection in the relative commutant of $P$. Relation \eqref{intro3} shows in particular that
\begin{equation}\label{intro5}
    L(\Sigma_n)\prec M_{n} \text{\; \; and \; \;} M_n\prec L(\Sigma_n).
\end{equation}
Finally, by building upon some results from \cite{DHI16}, we show that \eqref{intro4} and \eqref{intro5} imply that there exists a product decomposition $\Lambda=\Lambda_{\widehat n}\times\Lambda_n$ such that $M_{\widehat n}=L(\Lambda_{\widehat n})$ and $M_n=L(\Lambda_n)$, modulo unitary conjugacy and amplifications. 
The result follows now by induction.

{\bf Organization of the paper.} Besides the introduction there are four other sections in the paper. In Section 2 we review Popa's intertwining-by-bimodules techniques and some other tools. In Section 3 we recall some general properties of s-malleable deformations and some structural results for $L(\Gamma)$, when $\Gamma\in\mathcal A$.
Next, in Section 4 we review  Ioana's ultrapower technique and present an important corollary which, besides being of independent interest, will be used in the proof of Theorem \ref{A} . Finally, in Section 5 we prove our main results. 

{\bf Acknowledgments.} I am grateful to Adrian Ioana for helpful comments and to the anonymous referees for
their many remarks which greatly improved the exposition of the paper.

\section{Preliminaries}

\subsection{Terminology}
In this paper we consider {\it tracial von Neumann algebras} $(M,\tau)$, i.e. von Neumann algebras $M$ equipped with a faithful normal tracial state $\tau: M\to\mathbb C.$ This induces a norm on $M$ by the formula $\|x\|_2=\tau(x^*x)^{1/2},$ for all $x\in M$. We will always assume that $M$ is a {\it separable} von Neumann algebra, i.e. the $\|\cdot\|_2$-completion of $M$ denoted by $L^2(M)$ is separable as a Hilbert space.
We denote by $\mathcal Z(M)$ the {\it center} of $M$ and by $\mathcal U(M)$ its {\it unitary group}. For two von Neumann subalgebras $P_1,P_2\subset M$, we denote by $P_1\vee P_2=W^*(P_1\cup P_2)$ the von Neumann algebra generated by $P_1$ and $P_2$. 

All inclusions $P\subset M$ of von Neumann algebras are assumed unital. We denote by $E_{P}:M\to P$ the unique $\tau$-preserving {\it conditional expectation} from $M$ onto $P$, by $e_P:L^2(M)\to L^2(P)$ the orthogonal projection onto $L^2(P)$ and by $\langle M,e_P\rangle$ the Jones' basic construction of $P\subset M$. We also denote by $P'\cap M=\{x\in M|xy=yx, \text{ for all } y\in P\}$ the {\it relative commutant} of $P$ in $M$ and by $\mathcal N_{M}(P)=\{u\in\mathcal U(M)|uPu^*=P\}$ the {\it normalizer} of $P$ in $M$.  We say that $P$ is {\it regular} in $M$ if the von Neumann algebra generated by $\mathcal N_M(P)$ equals $M$. Next, we define the quasi-normalizer  $\mathcal {QN}_{M}(P)$ as the subgroup of all elements $x\in M$ for which there exist $x_1,...,x_n\in M$ such that $P x\subset \sum_i x_i P$ and $xP \subset \sum_i P x_i$ (see \cite[Definition 4.8]{Po99}). 

The {\it amplification} of a II$_1$ factor $(M,\tau)$ by a number $t>0$ is defined to be $M^t=p(\mathbb B(\ell^2(\mathbb Z))\bar\otimes M)p$, for a projection $p\in \mathbb B(\ell^2(\mathbb Z))\bar\otimes M$ satisfying $($Tr$\otimes\tau)(p)=t$. Here Tr denotes the usual trace on $\mathbb B(\ell^2(\mathbb Z))$. Since $M$ is a II$_1$ factor, $M^t$ is well defined. Note that if $M=P_1\bar\otimes P_2$, for some II$_1$ factors $P_1$ and $P_2$, then there exists a natural identification $M=P_1^t\bar\otimes P_2^{1/t}$, for every $t>0.$


For a countable group $\Gamma$ and for two subsets $S,T\subset \Gamma$, we denote 
by $C_S(T)=\{g\in S|gh=hg, \text{ for all } h\in T\}$ the {\it centralizer of $T$ in $S$}. Given a group inclusion $\Sigma<\Gamma$, 
the quasi-normalizer ${\rm QN}_\Gamma(\Sigma)$ is the group of all $g\in \Gamma$ for which exists a finite set $F\subset \Gamma$ such that $\Sigma g\subset F\Sigma$ and $g \Sigma\subset \Sigma F$; equivalently, $g\in {\rm QN}_\Gamma(\Sigma)$ if and only if $[\Sigma: g\Sigma g^{-1}\cap \Sigma]<\infty$ and $[\Sigma: g^{-1}\Sigma g\cap \Sigma]<\infty$. 

Finally, for a positive integer $n$, we denote by $\overline{1,n}$ the set $\{1,\dots, n\}$. For any subset $S\subset \overline{1,n}$ we denote its complement by $\widehat S=\overline{1,n}\setminus S$. If $S=\{i\},$ we will simply write $\hat i$ instead of $\widehat {\{i\}}$. Also, given
any product group $G=G_1\times \dots\times G_n$ and any tensor product $M_1\bar\otimes \dots\bar\otimes M_n$, we will denote their subproduct supported on $S$ by $G_S=\times_{i\in S}G_i$ and $M_S=\bar\otimes_{i\in S} M_i$, respectively.

\subsection {Intertwining-by-bimodules} We next recall from  \cite [Theorem 2.1 and Corollary 2.3]{Po03} the powerful {\it intertwining-by-bimodules} technique of S. Popa.

\begin {theorem}[\!\!\cite{Po03}]\label{corner} Let $(M,\tau)$ be a tracial von Neumann algebra and $P\subset pMp, Q\subset qMq$ be von Neumann subalgebras. Let $\mathcal U\subset\mathcal U(P)$ be a subgroup such that $\mathcal U''=P$.

Then the following are equivalent:

\begin{enumerate}

\item There exist projections $p_0\in P, q_0\in Q$, a $*$-homomorphism $\theta:p_0Pp_0\rightarrow q_0Qq_0$  and a non-zero partial isometry $v\in q_0Mp_0$ such that $\theta(x)v=vx$, for all $x\in p_0Pp_0$.

\item There is no sequence $(u_n)_n\subset\mathcal U$ satisfying $\|E_Q(xu_ny)\|_2\rightarrow 0$, for all $x,y\in M$.

\end{enumerate}
\end{theorem}

If one of the equivalent conditions of Theorem \ref{corner} holds true, we write $P\prec_{M}Q$, and say that {\it a corner of $P$ embeds into $Q$ inside $M$.}

\begin{notation} Throughout the paper we will use the following notation.
\begin{itemize}
\item If $Pp'\prec_{M}Q$ for any non-zero projection $p'\in P'\cap pMp$, then we write $P\prec^{s}_{M}Q$.
\item If $P\prec_{M}Qq'$ for any non-zero projection $q'\in Q'\cap qMq$, then we write $P\prec^{s'}_{M}Q$.
\item If both $P\prec^s_M Q$ and $P\prec^{s'}_M Q$ hold, then we write $P\prec_M^{s,s'} Q$.

\end{itemize}
\end{notation}




\subsection{From tensor decompositions to group product decompositions.}
In this subsection we review how some intertwining relations and a tensor product decomposition in a group von Neumann algebra can be used to deduce a direct decomposition of the underlying group, see Theorem \ref{Th:splits}. The result is actually a mild generalization of \cite[Theorem 6.1]{DHI16}. The proof follows the same idea, but for the convenience of the reader we include all the details. 

\begin{theorem}[\!\!\cite{DHI16}]\label{Th:splits}
Let $M=L(\Gamma)$ be the \text{II}$_1$ factor of an icc group $\Gamma$ and assume $M=P_1\bar\otimes P_2.$ Let $\Sigma_1,\Sigma_2<\Gamma$ be subgroups and $e\in L(\Sigma_1)'\cap M$ a non-zero projection such that
\begin{enumerate}
    \item $L(\Sigma_1)e\prec_M^s P_1$ and $P_1\prec_M L(\Sigma_1)e$,
    \item $L(\Sigma_2)\prec_M P_2$ and $P_2\prec_M L(\Sigma_2)$.
\end{enumerate}

Then there exist a decomposition $\Gamma=\Gamma_1\times \Gamma_2$, a decomposition $M=P_1^{t_1}\bar\otimes P_2^{t_2}$ for some $t_1,t_2>0$ with $t_1t_2=1$ and a unitary $u\in M$ such that $P_i^{t_i}=uL(\Gamma_i)u^*$, for any $i\in\overline{1,2}$. 
\end{theorem}

We first need  the following two lemmas inspired by \cite[Section 6]{DHI16}. Lemma \ref{L:joint} bellow can also be seen as an extension of \cite[Lemma 2.6]{Dr19}.

\begin{lemma}\label{QNbig}
Let $M=L(\Gamma)$ be the \text{II}$_1$ factor of an icc group $\Gamma$ and assume $M=P_1\bar\otimes P_2.$ Let $\Sigma<\Gamma$ be a subgroup and $e\in L(\Sigma)'\cap M$ a non-zero projection such that $L(\Sigma)e\prec_M^s P_1$ and $P_1\prec_M L(\Sigma)e$. 

Then $QN_\Gamma (\Sigma)<\Gamma$ has finite index.
\end{lemma}

{\it Proof.} 
We first use \cite[Lemma 2.4]{DHI16} and derive that there exists a non-zero projection $f\in \mathcal Z(L(\Sigma)'\cap M)$ satisfying $f\leq e$ and 
\begin{equation}\label{i11}
L(\Sigma)f\prec_M^s P_1 \text{ and }P_1\prec_M^{s'} L(\Sigma)f.    
\end{equation}
Note that $\mathcal{N}_{M}(L(\Sigma)'\cap M)'\cap M \subset L(\Sigma)'\cap M$. Hence, by passing to relative commutants in \eqref{i11}
it follows from \cite[Lemma 3.5]{Va08} and \cite[Lemma 2.4]{DHI16} that 
\begin{equation}\label{i2}
P_2\prec_M^{s'} (L(\Sigma)'\cap M)f \text{ and } (L(\Sigma)'\cap M)f\prec_M^s P_2.     
\end{equation}

By repeating the same argument, we deduce that $((L(\Sigma)'\cap M)'\cap M)f\prec_M^s P_1$. In combination with \eqref{i2} we derive that
$
\mathcal Z((L(\Sigma)'\cap M))f\prec_M^s P_1 \text{ and } \mathcal Z((L(\Sigma)'\cap M))f\prec_M^s P_2.
$
Therefore, by \cite[Lemma 2.8(2)]{DHI16} we further obtain that $\mathcal Z((L(\Sigma)'\cap M))f\prec_M^s \mathbb C1$, and hence, there exists a non-zero projection $f_1\in \mathcal Z((L(\Sigma)'\cap M))$ with $f_1\leq f$ such that $\mathcal Z((L(\Sigma)'\cap M))f_1=\mathbb C f_1.$

By letting $Q=L(\Sigma)\vee (L(\Sigma)'\cap M)$, it follows that $Q'\cap M\subset \mathcal Z (L(\Sigma)'\cap M)$, and therefore, $(Qf_1)'\cap f_1Mf_1=\mathbb Cf_1.$ 
Using relations \eqref{i11} and \eqref{i2} we have
$
P_1\prec_M Qf_1 \text{ and }P_2\prec_M Qf_1.
$
Therefore, we can apply \cite[Lemma 2.6]{Dr19} and derive that $M\prec_M Q$. Let $\Omega=\{g\in\Gamma|\; \mathcal O_{\Sigma}(g) \text{ is finite}\}$, where $\mathcal O_\Sigma(g)=\{hgh^{-1}|\; h\in \Sigma\}$ is the orbit of $g$ under conjugation by $\Sigma$. Note that $\Omega$ is normalized by $\Sigma$ and $L(\Sigma)'\cap M\subset L(\Omega)$. Since $\Omega\Sigma<{\rm QN}_\Gamma(\Sigma)$, it follows that $Q\subset L({\rm QN}_\Gamma(\Sigma))$. Altogether, it follows that $M\prec_M L({\rm QN}_\Gamma(\Sigma))$, which proves the lemma by using \cite[Lemma 2.5(1)]{DHI16}. 
\hfill$\blacksquare$


\begin{lemma}\label{L:joint}
Let $\Gamma\car B$ be a trace preserving action and denote $M=B\rtimes\Gamma$. Let $\Sigma_1<\Sigma_2<\Gamma$ be subgroups such that $L({QN}_{\Gamma}(\Sigma_1))'\cap M=\mathbb C1$.

Let $P_1, P_2\subset pMp$ be von Neuman subalgebras such that there exist commuting subalgebras $\tilde P_0,\tilde P_1, $ $\tilde P_2\subset M$ satisfying $P_1\subset \tilde P_1$, $P_2\subset \tilde P_2$ and $\tilde P_0\vee \tilde P_1\vee \tilde P_2=M$. 

If $P_i\prec_M B\rtimes\Sigma_i$, for any $i\in\overline{1,2}$,  then $P_1\vee P_2\prec_M B\rtimes\Sigma_2$.
\end{lemma}

{\it Proof.} Since $P_1\prec_M B\rtimes\Sigma_1$ and $P_2\prec_M B\rtimes\Sigma_2$, it follows that there exist projections $p_1\in P_1, p_2\in P_2$ and  non-zero elements $v_1 \in p_1M,v_2\in Mp_2$ such that
\begin{equation}\label{p1}
    (p_1P_1 p_1)_1v_1\subset v_1 (B\rtimes\Sigma_1)_1 \text{\; and \;} v_2 (p_2P_2 p_2)_1\subset (B\rtimes\Sigma_2)_1v_2. 
\end{equation}
We continue by showing that there exists $g\in QN_\Gamma(\Sigma_1)$ such that $v_1u_gv_2\neq 0$. If this is not the case, we derive that $v_1zv_2=0$, for all $z\in \mathcal U(L(QN_\Gamma(\Sigma_1)))$ and hence, $v_1zv_2v^*_2z^*=0$. By letting $a=\vee_{z\in \mathcal U(L(QN_\Gamma(\Sigma_1)))} zv_2v^*_2z$, it follows that $v_1a=0$. Since $a\in L(QN_\Gamma(\Sigma_1))'\cap M=\mathbb C1$, we derive that $v_1=0$ or $v_2=0$, contradiction.

Next, we consider $g\in QN_\Gamma(\Sigma_1)$ such that $v_1u_gv_2\neq 0$. This implies that there exist some elements $g_1,\dots,g_n\in \Gamma$ such that 
$(B\rtimes\Sigma_1)_1u_g\subset \sum_{j=1}^n u_{g_j} (B\rtimes\Sigma_1)_1.$ Using \eqref{p1}, we derive that
\begin{equation}\label{p2}
    \mathcal U(p_1P_1p_1 ) (v_1 u_g v_2) \mathcal U(p_2P_2p_2)\subset \sum_{i=1}^n v_1 u_{g_j} (B\rtimes\Sigma_2)_1v_2.
\end{equation}
Now, we assume by contradiction that $p_1P_1p_1\vee p_2P_2p_2\nprec_M B\rtimes\Sigma_2$. Hence, there exist two sequence of unitaries $(a_n)_n\subset \mathcal U(p_1P_1p_1)$ and $(b_n)_n\subset \mathcal U(p_2P_2p_2)$ such that $\|E_{B\rtimes\Sigma_2}(xa_n b_n y)\|_2\to 0$, for all $x,y\in M$. Since $P_1\subset \tilde P_1$, $P_2\subset \tilde P_2$ and $\tilde P_0\vee \tilde P_1\vee \tilde P_2=M$, we further derive that $\|E_{B\rtimes\Sigma_2}(xa_n (v_1zv_2) b_n y)\|_2\to 0$, for all $x,y\in M$. 

For any subset $F\subset\Gamma$, we denote by $P_F$ the orthogonal projection onto the closed linear span of $\{Bu_g|g\in F\}$ and notice that $\|P_{G\Sigma_2 H}(x)\|_2^2\leq \sum_{g\in G,h\in H} \|E_{B\rtimes\Sigma_2}(u_g^* xu_h ^*)\|_2^2$, for all finite subsets $G,H\subset \Gamma$  and $x\in M$.
Hence, the previous paragraph implies that $\|P_{G\Sigma_2 H}(a_n (v_1zv_2) b_n)\|_2\to 0$, for all finite subsets $G,H\subset \Gamma.$ By using Kaplansky's density theorem, this contradicts \eqref{p2} and therefore ends the proof of the lemma.
\hfill$\blacksquare$



{\it Proof of Theorem \ref{Th:splits}.}
From assumption we have $L(\Sigma_1)\prec_M L(\Sigma_2)'\cap M$. By applying \cite[Theorem 6.2]{DHI16} we obtain that there exist finite index subgroups $\Theta_1<k\Sigma_1k^{-1}$ and $\Theta_2<\Sigma_2$ for some $k\in\Gamma$ such that the commutator subgroup $[\Theta_1,\Theta_2]$ is finite. Denote by $\Theta<\Gamma$ the subgroup generated by $\Theta_1$ and $\Theta_2$.

By Lemma \ref{QNbig} we have $[\Gamma: {\rm QN}_\Gamma(\Sigma_1)]<\infty$.
Since $\Gamma$ is icc, we derive that $L({\rm QN}_\Gamma(\Theta_1))'\cap M= \mathbb C1$. Therefore, by Lemma \ref{L:joint} we derive that $M\prec_M L(\Theta)$, which shows that $[\Gamma:\Theta]<\infty$. In particular, $\Theta$ is icc. Since $[\Theta_1,\Theta_2]$ is a finite normal subgroup of $\Theta$, it must be trivial. Therefore, $\Theta_1$ and $\Theta_2$ are commuting subgroups of $\Gamma$ such that $\Theta_1\cap\Theta_2=1$ and $[\Gamma:\Theta_1\Theta_2]<\infty$.

Finally, by applying verbatim the second part of the proof of \cite[Theorem 6.1]{DHI16} (or the second part of the proof of \cite[Theorem 4.14]{CdSS15}), we obtain the conclusion of the theorem.
\hfill$\blacksquare$

\subsection{An augmentation technique for intertwining}
An essential tool for the proof of our main result, Theorem \ref{A}, is the augmentation technique developed
in \cite{CD-AD20}; this was also crucial for the recent work \cite{CD-AD21}. One place where the augmentation technique is used in our paper is in the following result.

\begin{proposition}\label{augm}
Let $\Gamma$ be a countable icc group and let $M=L(\Gamma)$. Let $\Sigma,\Theta<\Gamma$ be subgroups and $M=P\bar\otimes Q$ a tensor product decomposition. Let $G_0<G$ be countable groups and let $G\car A$ be a trace preserving action such that $Q=A\rtimes G$. Assume that the following two conditions hold:
\begin{itemize}
    \item $P\prec_M L(\Sigma)$ and $L(\Sigma)\prec_M P$.
    
    \item $Q\prec_M L(\Theta)$ and $L(\Theta)\prec_M P\bar\otimes (A\rtimes G_0)$.
\end{itemize}

Then $[G:G_0]<\infty$.
\end{proposition}

{\it Proof.}
Following the augmentation technique from \cite[Section 3]{CD-AD20}, we consider a Bernoulli action $\Gamma\car D$  with abelian base. Denote $\mathcal M=D\rtimes\Gamma$ and let $\Psi:\mathcal M\to \mathcal M\bar\otimes M$ be the $*$-homomorphism given by $\Psi(du_g)=du_g\otimes u_g$, for all $d\in D$ and $g\in\Gamma$. 

From \cite[Lemma 2.4(2)]{DHI16} we get that $P\prec_M^s L(\Sigma)$, and hence, $P\prec_{\mathcal M}^{s,s'} D\rtimes\Sigma$. By using \cite[Remark 2.2]{DHI16} and \cite[Lemma 2.3]{Dr19} we further obtain that $\Psi(P)\prec^{s,s'}_{\mathcal M \bar\otimes M}\mathcal M\bar\otimes L(\Sigma)$. Using $L(\Sigma)\prec_M P$ and \cite[Lemma 2.4(2)]{Dr19}, we get $\Psi(P)\prec^{s}_{\mathcal M \bar\otimes M}\mathcal M\bar\otimes P$. In the same way, we derive that $\Psi(Q)\prec^{s}_{\mathcal M \bar\otimes M}\mathcal M\bar\otimes P\bar\otimes (A\rtimes G_0)$. Using \cite[Lemma 2.3]{BV12} we get that $\Psi(M)\prec_{\mathcal M \bar\otimes M}\mathcal M\bar\otimes P\bar\otimes  (A\rtimes G_0)$. Since $\Psi(x)=x\otimes 1$, for all $x\in D$, we deduce that $\Psi(\mathcal M)\prec_{\mathcal M \bar\otimes M}\mathcal M\bar\otimes P\bar\otimes  (A\rtimes G_0)$. By \cite[Lemma 10.2]{IPV10} we obtain that $M\prec_{\mathcal M} P\bar\otimes  (A\rtimes G_0)$. Since $\mathcal {QN}_{\mathcal M}^{(1)}(M)=M$, we derive that $M\prec_{ M} P\bar\otimes  (A\rtimes G_0)$, which implies from \cite[Lemma 2.5(1)]{DHI16} that $[G:G_0]<\infty$.
\hfill$\blacksquare$

\subsection{Relative amenability and weak containment of bimodules}

A tracial von Neumann algebra $(M,\tau)$ is {\it amenable} if there exists a positive linear functional $\Phi:\mathbb B(L^2(M))\to\mathbb C$ such that $\Phi_{|M}=\tau$ and $\Phi$ is $M$-{\it central}, meaning $\Phi(xT)=\Phi(Tx),$ for all $x\in M$ and $T\in \mathbb B(L^2(M))$. By Connes' breakthrough classification of amenable factors \cite{Co76}, it follows that $M$ is amenable if and only if $M$ is approximately finite dimensional.
 
Next, we recall the notion of relative amenability introduced by 
Ozawa and Popa in \cite{OP07}. Let $(M,\tau)$ be a tracial von Neumann algebra. Let $p\in M$ be a projection and $P\subset pMp,Q\subset M$ be von Neumann subalgebras. Following \cite[Definition 2.2]{OP07}, we say that $P$ is {\it amenable relative to $Q$ inside $M$} if there exists a positive linear functional $\Phi:p\langle M,e_Q\rangle p\to\mathbb C$ such that $\Phi_{|pMp}=\tau$ and $\Phi$ is $P$-central. 
Note that $P$ is amenable relative to $\mathbb C$ inside $M$ if and only if $P$ is amenable.
We say that $P$ is {\it strongly non-amenable relative to} $Q$ if $Pp'$ is non-amenable relative to $Q$ for any non-zero projection $p'\in P'\cap pMp$.  

Let $M, N$ be tracial von Neumann algebras. An $M$-$N$ bimodule $_M \mathcal H_N$ is a Hilbert space $\mathcal H$ equipped with a $*$-homomorphism $\pi_{\mathcal H}: M\odot N^{\text{op}}\to\mathbb B(\mathcal H)$ that is normal on $M$ and $N^{\text{op}}$, where $M\odot N^{\text{op}}$ is the algebraic tensor product between $M$ and the oposite von Neumann algebra $N^{\text{op}}$ of $N$. For two $M$-$N$ bimodules $_M \mathcal H_N$ and $_M \mathcal K_N$, we say that $_M \mathcal H_N$ is weakly contained in $_M \mathcal K_N$ if $\|\pi_{\mathcal H}(x)\|\leq \|\pi_{\mathcal K}(x)\|$, for any $x\in M\odot N^{\text{op}}$. Examples of bimodules include the trivial bimodule $_ML^2(M)_M$ and the coarse bimodule $_ML^2(M)\otimes L^2(N)_N$.

Finally, if $Q\subset M$ are tracial von Neumann algebras, then $_M L^2(\langle M,e_{Q} \rangle)_M\cong$$_M( L^2(M)\otimes_Q L^2(M))_M$. Also, a subalgebra $P\subset pMp$ is amenable relative to $Q$ if and only if $_P L^2(pM)_M$ is weakly contained in $_P L^2(\langle M,e_{Q} \rangle)_M$.

\subsection{Mixing inclusion of von Neumann algebras}

An important tool in the intertwining-by-bimodules techniques is the mixing notion relative to a subalgebra (see \cite[Definition 2.9]{Po05}, \cite[Definition 2.3]{PS09} and \cite[Definition A.4.2]{Bo14}).

\begin{definition}\label{def:mixing}
Let $A\subset M\subset \tilde M$ be an inclusion of tracial von Neumann algebras. We say that $M\subset \tilde M$ is {\it mixing relative} to $A$ if for any sequence $(u_n)_n\subset (M)_1$ satisfying $\|E_{A}(xu_n y)\|_2\to 0$, for all $x,y\in M$, we have $\|E_{M}(\tilde x u_n \tilde y)\|_2\to 0$, for all $\tilde x,\tilde y\in \tilde M\ominus M$.
\end{definition}

In Definition \ref{def:mixing}, if $A=\mathbb C1$ we simply say that $M\subset\tilde M$ is {\it mixing}.
We record the following well known lemmas and include the proof of the first one only for the convenience of the reader; the second one can be proven in a similar way. 

\begin{lemma}\label{mixingbernoulli}
Let $\Sigma<\Gamma$ be countable groups and denote $I=\Gamma/\Sigma$. Let $A_0\subset B_0$ be tracial von Neumann algebras and denote $M=A_0^{I}\rtimes\Gamma $ and $\tilde M=B_0^{I}\rtimes\Gamma $. Note that $M\subset\tilde M$.

Then $M\subset \tilde M$ is mixing relative to $A_0^{I}\rtimes\Sigma.$ 
\end{lemma}

{\it Proof.} Let $(w_n)_n\subset\mathcal U(M)$ be a sequence of unitaries such that $\|E_{A_0^{I}\rtimes\Sigma}(xw_ny)\|_2\to 0$, for all $x,y\in M$. We have to show that $\|E_{M}(\tilde x w_n \tilde y)\|_2\to 0$, for all $\tilde x,\tilde y\in \tilde M\ominus M$. By Kaplansky's density theorem, it is enough to assume $\tilde x=b\in (B_0\ominus A_0)^{h\Sigma}$ and $\tilde y=c\in (B_0\ominus A_0)^{k\Sigma}$ for some $h,k\in I$. If we let $w_n=\sum_{g\in\Gamma} w_n^gu_g \in A_0^I\rtimes\Gamma$, note that $E_M(bw_n^g\sigma_g(c))=0$ if $g\notin  h\Sigma k^{-1}$ and hence 
\[
\|E_M(\tilde xw_n\tilde y)\|_2\leq \|b\| \|c\| \|\sum_{g\in h\Sigma k^{-1}}w_n^gu_g\|_2 = \|b\| \|c\| \| E_{A_0^I\rtimes\Sigma}(u_{h^{-1}}w_nu_k)\|_2 \to 0,
\]
by the assumption.
\hfill$\blacksquare$

\begin{lemma}\label{mixingamalgam}
Let $\tilde M=M_1*_AM_2$ be an amalgamated free product of tracial von Neumann algebras. Then $M_1\subset\tilde M$ is mixing relative to $A$.
\end{lemma}

\section{Two classes of von Neumann algebras that admit s-malleable deformations}

\subsection{Malleable deformations.}
In \cite{Po01, Po03} Popa introduced the notion of an s-malleable deformation of a von Neumann algebra. In the framework of his powerful deformation/rigidity techniques, this notion has led to a remarkable progress in the theory of von Neumann algebras, see the surveys \cite{Po07,Va10a, Io12b,Io17}. See also \cite{dSHHS20} for a comprehensive overview on s-malleable deformations and for recent developments.  

\begin{definition}\label{def:malleable}
Let $(M,\tau)$ be a tracial von Neumannn algebra. A pair $(\tilde M, (\alpha_t)_{t\in\mathbb R})$ is called an {\it s-malleable deformation} of $M$ if the following conditions hold:
\begin{itemize}
    \item $(\tilde M,\tilde\tau)$ is a tracial von Neumann algebra such that $M\subset\tilde M$ and $\tau=\tilde\tau_{|M}.$ 
    \item $(\alpha_t)_{t\in\mathbb R}\subset {\text Aut}(\tilde M,\tilde\tau)$ is a $1$-parameter group with $\lim_{t\to 0}\|\alpha_t(x)-x\|_2=0$, for any $x\in\tilde M.$ 
    \item There exists $\beta\in {\text Aut}(\tilde M,\tilde\tau)$ that satisfies $\beta_{|M}=\text{Id}_M$, $\beta^2=\text{Id}_{\tilde M}$ and $\beta\alpha_t=\alpha_{-t}\beta$, for any $t\in\mathbb R$.
    
          \item $\alpha_t$ does not converge uniformly to the identity on $(M)_1$ as $t\to 0$.

\end{itemize}

\end{definition}

\begin{theorem}[\!\!\cite{dSHHS20}]\label{Th:deform}
Let $(\tilde M, (\alpha_t)_{t\in\mathbb R})$ be an s-malleable deformation of a tracial von Neumann algebra $M$. Let $A\subset M$ and $Q\subset qMq$ be some von Neumann subalgebras such that $M\subset \tilde M$ is mixing relative to $A$, $\alpha_t\to \text{id}$ uniformly on $(Q)_1$ and $Q\nprec_M A$. Then the following hold:
\begin{enumerate}
    \item  $\alpha_t\to \text{id}$ uniformly on  $(Q\vee (Q'\cap qMq))_1$. 
    \item Let $Q=Q_1\subset Q_2\subset\dots\subset qMq$ be an ascending sequence of von Neumann subalgebras with $\alpha_t\to \text{id}$ uniformly on $(Q_i)_1$, for all $i\ge 1$. Then $\alpha_t\to \text{id}$ uniformly on  $(\bigvee_{i\ge 1} Q_i)_1$.
    
\end{enumerate}
\end{theorem}

{\it Proof.} Part (1) follows from \cite[Corollary 6.7(ii)]{dSHHS20}. For part (2), by using \cite[Lemma 6.3]{dSHHS20} it follows that $Q'\cap q\tilde M q\subset qMq$ and the conclusion follows from \cite[Theorem 3.5]{dSHHS20}.
\hfill$\blacksquare$

We will also need the following result  from \cite[Proposition 5.6]{dSHHS20}.

\begin{proposition}[\!\!\cite{dSHHS20}]\label{P:central}
Let $(\tilde M, (\alpha_t)_{t\in\mathbb R})$ be an s-malleable deformation of a tracial von Neumann algebra $M$. Let $Q\subset qMq$ be a von Neumann subalgebra and let $q_0\in Q$ be a non-zero projection such that $\alpha_t\to \text{id}$ uniformly on $(q_0Qq_0)_1$. 

Then $\alpha_t\to \text{id}$ uniformly on $(Qz)_1$, where $z$ is the central support of $q_0$ in $Q$. 
\end{proposition}

\subsection{\bf Class $\mathscr M$}\label{class M}
We say that a non-amenable II$_1$ factor $M$ is in Class $\mathscr M$ if there exists an s-malleable deformation $(\tilde M, (\alpha_t)_{t\in\mathbb R})$ of $M$ and an amenable subalgebra $A\subset M$ satisfying:
\begin{enumerate}
        \item The inclusion $M\subset\tilde M$ is mixing relative to $A$.
    
    \item For any tracial von Neumann algebra $N$ and for any subalgebra $P\subset p(M\bar\otimes N)p$ such that $P'\cap p(M\bar\otimes N)p$ is strongly non-amenable relative to $1\otimes N$, it follows that:
    \begin{enumerate}[(i)]
        \item $\alpha_t\otimes \text{id}\to \text{id} \text{ uniformly on }(P)_1.$
        \item If $P\prec_{M\otimes N} A\otimes N,$ then $P\prec_{M\otimes N} 1\otimes N.$
    \end{enumerate}
    
    
\end{enumerate}

    
    



While this class of II$_1$ factors seems somewhat technical, it actually contains all group von Neumann algebras $L(\Gamma)$ with $\Gamma\in\mathcal A$ and all non-trivial tracial free products $M_1*M_2$, see Proposition \ref{examples} bellow and its proof. Note also that if $A=\mathbb C 1$, then condition (2) is simply reflecting Popa's spectral gap principle.

\begin{proposition}\label{examples}
If $\Gamma\in\mathcal A$, then $L(\Gamma)\in\mathscr M$.



\end{proposition}

{\it Proof.} 
If $\Gamma\in\mathcal A_1$ we recall that Sinclair constructed in \cite[Section 3]{Si10} an s-malleable deformation $(\tilde M, (\alpha_t)_{t\in\mathbb R})$ in the sense of Definition \ref{def:malleable}; see also \cite[Section 3.1]{Va10b} and \cite[Section 2]{Io11}. We will prove that $M\in\mathscr M$ with $A=\mathbb C1$. Note first that $M\subset \tilde M$ is mixing since $\pi$ is a mixing representation. Next, Lemma \ref{L:spectralgapcohomology} bellow shows that condition (2) of Class $\mathscr{M}$ is satisfied.

Next, if $\Gamma=\Gamma_1*_\Sigma\Gamma_2\in\mathcal A_2$, we recall that \cite[Section 2.2]{IPP05} shows that $M=L(\Gamma)$ admits an s-malleable deformation $(\tilde M, (\alpha_t)_{t\in\mathbb R})$ in the sense of Definition \ref{def:malleable} with $\tilde M=M*_{L(\Sigma)}(L(\Sigma)\bar\otimes L(\mathbb F_2))$. It follows that $M\in\mathscr{M}$ with $A= L(\Sigma).$ Indeed, we note first that $M\subset \tilde M$ is mixing relative to $A$ by Lemma \ref{mixingamalgam}.
Next, since we have the decomposition $M\bar\otimes N= (L(\Gamma_1)\bar\otimes N)*_{L(\Sigma)\bar\otimes N} (L(\Gamma_2)\bar\otimes N)$, condition (2.i) follows from \cite[Lemma 6.5]{Io12a}. To show condition (2.ii), let $N$ be a tracial von Neumann algebra and $P\subset p(M\bar\otimes N)p$  a subalgebra such that $P'\cap p(M\bar\otimes N)p$ is strongly non-amenable relative to $1\otimes N$. If $P\prec_{M\otimes N} A\otimes N$ and $P\nprec_{M\otimes N} 1\otimes N,$ then we derive (see, e.g., \cite[Proposition 3.7]{Dr17}) that $P'\cap p(M\bar\otimes N)p\prec_{M\bar\otimes N} A\bar\otimes N$. It implies by \cite[Lemma 2.6(3)]{DHI16} that $P'\cap p(M\bar\otimes N)p$ is not strongly non-amenable relative to $1\otimes N$, contradiction. 

Finally, if $\Gamma\in\mathcal A_3$, we recall that Ioana constructed in \cite[Section 2]{Io06} an s-malleable deformation $(\tilde M, (\alpha_t)_{t\in\mathbb R})$ in the sense of Definition \ref{def:malleable}, see Remark \ref{remarkbernoullimalleable}(1).
Next, we note that $M\subset \tilde M$ is mixing relative to $A:=L(\Sigma)^{G/H}\rtimes H$ by Lemma \ref{mixingbernoulli}. Next, parts (i) and (ii) from condition (2) of Class $\mathscr M$ follow from \cite[Corollary 4.3]{IPV10} and its proof. 
\hfill$\blacksquare$

The following lemma is a standard application of Popa's spectral gap rigidity principle \cite{Po06b} and it is essentially contained in the proof of \cite[Theorem 6.4]{Ho15} (see also \cite[Lemma 2.2]{Io11}). For completeness, we include a proof.

\begin{lemma}\label{L:spectralgapcohomology}
Let $\Gamma$ be a countable non-amenable group that admits an unbounded cocycle for some mixing representation $\pi:\Gamma\to \mathcal O( H_{\mathbb R})$ such that $\pi$ is weakly contained in the left regular representation of $\Gamma$. 
Let $\Gamma\car N$ be a trace preserving action  and denote $M=N\rtimes\Gamma$. 

If $P\subset pMp$ is a von Neumann subalgebra that is strongly non-amenable relative to $N$, then $\alpha_t\to id$ uniformly on $(P'\cap pMp)_1.$
\end{lemma}

{\it Proof.} Since $\pi$ is contained in the left regular representation, it follows by \cite[Lemma 3.5]{Va10b} that the
 $M$-$M$ bimodule $L^2(\tilde M \ominus M)$ is weakly contained in the $M$-$M$ bimodule $L^2(M)\otimes_N L^2(M)$.  
Since $P$ is strongly non-amenable relative to $N$, it follows that for any non-zero projection $z\in \mathcal Z(P'\cap pMp)$ we have that $_M L^2(M)_{Pz}$ is not weakly contained in $_M L^2(M)\otimes_N L^2(M)_{Pz}$. Hence, for any non-zero projection $z\in \mathcal Z(P'\cap pMp)$ we derive that $_M L^2(M)_{Pz}$ is not weakly contained in $_M L^2(\tilde M \ominus M)_{Pz}$. 

Let $\epsilon>0$. Therefore, by \cite[Lemma 2.3]{IPV10} we obtain that there exist $a_1,\dots,a_n\in P$ and $\delta>0$ such that if $x\in (p\tilde M p)_1$ satisfies $\|xa_i-a_i x\|_2\leq \delta$, for any $i\in\overline{1,n},$ then $\|x-E_{M}(x)\|_2\leq \epsilon$. We choose $t_0>0$ such that $\|\alpha_t(a_i)-a_i\|_2\leq\delta/2$, for all $|t|\leq t_0$ and $i\in\overline{1,n}$. Let $x\in (P'\cap pMp)_1$ and $t\in\mathbb R$ such that $|t|\leq t_0$. By using the triangle inequality we derive that for all $i\in\overline{1,n}$ and $|t|\leq t_0$ we have
\[
\|a_i\alpha_t(x)-\alpha_t(x)a_i \|_2= \|\alpha_{-t}(a_i)x-x\alpha_{-t}(a_i) \|_2 \leq 2\|\alpha_{-t}(a_i)-a_i\|_2\leq\delta.
\]
As a consequence, we deduce that $\|\alpha_t(x)-E_M(\alpha_t(x))\|_2\leq \epsilon$, for any $|t|\leq t_0$. Using Popa's transversality property \cite[Lemma 3.1]{Va10b} we obtain that 
\[
\|\alpha_t(x)-x\|_2\leq \sqrt{2}\epsilon, \text{ for all } x\in (P'\cap pMp)_1 \text{ and } |t|\leq t_0.
\]
This concludes the proof. \hfill$\blacksquare$


\begin{remark}\label{inclusion}
The proof of Lemma \ref{L:spectralgapcohomology} shows that the class $\mathscr M_0$ defined in Remark \ref{newrigidity} is contained in $\mathscr{M}$ and that any non-amenable tracial free product $M_1*M_2$ belongs to $\mathscr M_0$ (see also \cite[Lemma 2.10]{Io12a}).
\end{remark}

We end this subsection by showing that any von Neumann algebra $M\in\mathscr M$ does not have  property Gamma, i.e. for any uniformly bounded sequence $(x_n)_n\subset M$ with $\|x_ny-yx_n\|_2\to 0$, for any $y\in M$, must satisfy $\|x_n-\tau(x_n)\|_2\to 0.$

\begin{lemma}
Let $M$ be a tracial von Neumann algebra that belongs to Class $\mathscr M$. Then $M$ does not have property Gamma.
\end{lemma}

{\it Proof.} We assume the contrary. By using \cite[Theorem 3.1]{HU15} it follows that there exists a decreasing sequence of diffuse abelian von Neumann subalgebras $A_n\subset M$ with $n\ge 1$ such that $M=\bigvee_{n\ge 1} (A_n'\cap M$). 

Let $(\tilde M, (\alpha_t)_{t\in\mathbb R})$ be an s-malleable deformation of $M$ and $A\subset M$ a subalgebra as given by the definition of Class $\mathscr{M}$. Since $M$ is non-amenable, it follows that there exists $n\ge 1$ such that $A_n'\cap M$ is non-amenable. Let $z\in \mathcal Z(A_n'\cap M)$ such that $(A_n'\cap M)z$ is strongly non-amenable relative to $\mathbb C 1$. By using the fact that $M$ belongs to Class $\mathscr{M}$, it follows that $\alpha_t\to\text{id}$ uniformly on $(A_nz)_1$. In particular,
\begin{equation}\label{s1}
\alpha_t\to\text{id} \text{ uniformly on }(A_mz)_1, \text{ for any }m\ge n.     
\end{equation}
Note here that $z\in A_m'\cap M$, for any $m\ge n$.
If $A_mz\prec_M A$, condition (2.ii) of Class $\mathscr{M}$ shows that there exists a non-zero projection $z_1\in (A_m'\cap M)'\cap M\subset (A_n'\cap M)'\cap M$ with $z_1\leq z$ such that $(A_m'\cap M)z_1$ is amenable. This proves that $(A_n'\cap M)z_1$ is amenable,  contradiction.

Hence, $A_mz\nprec_M A$. Since $M\subset\tilde M$ is mixing relative to $A$, 
it follows from \eqref{s1} and Theorem \ref{Th:deform}(1) that $\alpha_t\to\text{id}$ uniformly on $(z(A_m'\cap M)z)_1$, for any $m\ge n$. Since $z(A_m'\cap M)z\nprec_M A$ and $zMz= \bigvee_{m\ge n} z(A_m'\cap M)z$, we apply Theorem \ref{Th:deform}(2) and derive that $\alpha_t\to\text{id}$ uniformly on $(zMz)_1$. Since $M$ is a factor, we apply Proposition \ref{P:central} and deduce that $\alpha_t\to\text{id}$ uniformly on $(M)_1$, contradiction. Therefore, $M$ does not have property Gamma.
\hfill$\blacksquare$

\subsection {Class $\mathscr{M}_{wr}$.}
We say that a von Neumann algebra $M$ is in Class $\mathscr{M}_{wr}$ if there exists a decomposition $M=B^I_0\rtimes\Lambda$ satisfying the following properties:
\begin{itemize}
    \item $B_0$ is a tracial amenable von Neumann algebra and $\Lambda$ is a non-amenable group.
    \item There exists $k\ge 1$ such that Stab$\, J$ is finite whenever $J\subset I$ with $|J|\ge k$.
\end{itemize}

\begin{remark}\label{remarkbernoullimalleable}
We record the following properties of von Neumann algebras that belong to $\mathscr M_{wr}$.
\begin{enumerate}
    \item Throughout the proofs of the main results, we will use the fact that
    any $M\in\mathscr M_{wr}$ admits an s-malleable deformation $(\tilde M, (\alpha_t)_{t\in\mathbb R})$ in the sense of Definition \ref{def:malleable} by using the free product deformation, see \cite[Section 2]{Io06}. To recall this construction, we define a self-adjoint unitary $h\in L(\mathbb Z)$ with spectrum $[-\pi,\pi]$ such that exp$(ih)$ equals the canonical generating unitary $u\in L(\mathbb Z)$. For any $t\in\mathbb R$, define $u_t=$exp$(ith)\in L(\mathbb Z)$. We let $\tilde M=L(B_0*L(\mathbb Z))^{I}\rtimes \Lambda\supset M$ and $\alpha_t=\otimes_{i\in I}\text{Ad}(u_t)\in\text{Aut}(\tilde M)$.
    \item If $M\in\mathscr M_{wr}$, then $M$ is a II$_1$ factor without property Gamma (see, e.g., \cite[Proposition 4.3]{Dr20}).
\end{enumerate}
\end{remark}

\section{From commuting subalgebras to commuting subgroups}
One of the crucial ingredients of the proof of Theorem \ref{A} is an ultrapower technique due to Adrian Ioana \cite{Io11}, which we recall in the following form. This result is essentially contained in the proof of \cite[Theorem 3.1]{Io11} (see also \cite[Theorem 3.3]{CdSS15}). The statement that we will use is a particular case of \cite[Theorem 4.1]{DHI16}.

\begin{theorem}[\!\!{\cite{Io11}}]\label{Th:ultrapower}
Let $\Gamma$ be a countable icc group and denote by $M=L(\Gamma)$. Let $\Delta:M\to M\bar\otimes M$ be the $*$-homomorphism given by $\Delta (u_g)=u_g\otimes u_g$, for all $g\in\Gamma.$ Let $B,Q\subset M$ be von Neumann subalgebras such that $\Delta(B)\prec_{M\bar\otimes M} M\bar\otimes Q$.

Then there exists a decreasing sequence of subgroups $\Sigma_k<\Gamma$ such that $B\prec_M L(\Sigma_k)$, for every $k\ge 1$, and $Q'\cap M\prec_M L(\cup_{k\ge 1} C_\Gamma(\Sigma_k)).$
\end{theorem}

Note that the ultrapower technique \cite{Io11} has been of crucial use in several other works in order to obtain structural and rigidity results for certain classes of group and group measure space von Neumann algebras \cite{CdSS15, KV15, DHI16, CI17, CU18, Dr19, CDK19, CDHK20, CD-AD20,CD-AD21}.

By combining Theorem \ref{Th:ultrapower} with a recent characterization of von Neumann algebras that do not have property Gamma \cite{BMO19,IM19} we obtain the following useful consequence.

\begin{theorem}\label{Cor:ultrapower}
We consider the context of Theorem \ref{Th:ultrapower}. In addition, assume that $M=P\bar\otimes Q$ where $P$ does not have property Gamma. 

Then there exists a subgroup $\Sigma<\Gamma$ with non-amenable centralizer  $C_\Gamma(\Sigma)$ such that $B\prec_M L(\Sigma)$ and $L(\Sigma)\prec_M Q.$ 
\end{theorem}

{\it Proof.} The proof is inspired by \cite[Lemma 5.2]{IM19}. We start the proof by recalling two general facts from \cite{BMO19}. First, \cite[Proposition 3.2]{BMO19} is stating that if a tracial von Neumann algebra $N$ does not have property Gamma, then any $N$-$N$ bimodule $\mathcal K$ that is weakly equivalent to $L^2(N)$ (i.e. $_N \mathcal K_N$ is weakly contained in $L^2(N)$ and $L^2(N)$ is weakly contained in $_N\mathcal K_N$)  must contain $L^2(N)$. 

Second, by using \cite[Lemma 3.4]{BMO19} we derive that for any subalgebra $Q_0\subset M$ we have
\begin{equation}\label{zx1}
L^2(\langle M,e_{Q_0'\cap M}\rangle) \text{ is weakly contained in }L^2(M) \text{ as } M\text{-}M \text{ bimodules}.     
\end{equation}
Next, we consider a decreasing sequence of subgroups $\Sigma_k<\Gamma$ as in the conclusion of Theorem \ref{Th:ultrapower}. Since $P$ is non-amenable, we can assume without loss of generality that $C_\Gamma(\Sigma_k)$ is non-amenable for any $k\ge 1.$
Denote the $M$-$M$ bimodule $\mathcal H=\oplus_{k\ge 1} L^2(\langle M, e_{L(\Sigma_k)'\cap M} \rangle)$ and notice that \eqref{zx1} implies that
\begin{equation}\label{zx2}
\mathcal H  \text{ is weakly contained in } L^2(M) \text{ as }M\text{-}M \text{ bimodules.}    
\end{equation}
On the other hand, since $P\prec_M L(\cup_{k\ge 1} C_\Gamma(\Sigma_k))$, it follows by Lemma 2.4(2) and Lemma 2.6(3) from \cite{DHI16} that $P$ is amenable relative to $\bigvee_{k\ge 1}(L(\Sigma_k)'\cap M)$ inside $M$. Hence, $L^2(M)$ is weakly contained in $L^2(\langle M,e_{\bigvee_{k\ge 1}(L(\Sigma_k)'\cap M)}\rangle)$ as $M$-$P$ bimodules. By using the moreover part of \cite[Proposition 2.5]{IM19}, we deduce that $L^2(\langle M,e_{\bigvee_{k\ge 1}(L(\Sigma_k)'\cap M)}\rangle)$ is weakly contained in $\mathcal H$ as $M$-$M$ bimodules. Therefore,
\begin{equation}\label{zx3}
    L^2(M) \text{ is weakly contained in } \mathcal H \text{ as }M\text{-}P \text{ bimodules.} 
\end{equation}
Since $M=P\bar\otimes Q$, it follows that $L^2(M)$ is a multiple of $L^2(P)$ as $P$-$P$ bimodules. In combination with \eqref{zx2} and \eqref{zx3}, we deduce that 
$L^2(P) \text{ is weakly equivalent to } \mathcal H  \text{ as }P\text{-}P \text{ bimodules.}$ Finally, since $P$ does not have property Gamma, we use the first paragraph of the proof and derive that $L^2(P)$ is contained in $_P \mathcal  H_P.$ Hence, there exists $k_0\ge 1$ such that $L^2(P)$ is contained in $_P L^2(M,e_{L(\Sigma_{k_0})'\cap M})_P$. This implies that $P\prec_M L(\Sigma_{k_0})'\cap M$. By passing to relative commutants and using \cite[Lemma 3.5]{Va08}, the conclusion of the corollary is obtained by taking $\Sigma=\Sigma_{k_0}.$
\hfill$\blacksquare$

As an application of Theorem \ref{Cor:ultrapower}, we can use
the augmentation technique from \cite[Section 3]{CD-AD20} and derive the following result which is an important ingredient of the proof of Theorem \ref{A}.

\begin{theorem}\label{Th:generalonebyone}
Let $\Gamma$ be a countable icc group such that $L(\Gamma)$ does not have property Gamma. Denote $M=L(\Gamma)^{1/t}$ for some $t>0$ and let $\Delta:M^t\to M^t\bar\otimes M^t$ be the $*$-homomorphism given by $\Delta(u_g)=u_g\otimes u_g$, for any $g\in\Gamma$. 

Assume  $M=M_1\bar\otimes\dots \bar\otimes M_n$ is the tensor product of $n\ge 1$ II$_1$ factors with the property that for any $i\in\overline{1,n}$, there exists $f(i)\in\overline{1,n}$ such that $\Delta(M^t_{\widehat i})\prec^s_{M^t\bar\otimes M^t} M^t\bar\otimes M^t_{\widehat {f(i)}}$.

Then $\Delta(M^t_{I})\prec^s_{M^t\bar\otimes M^t} M^t\bar\otimes M^t_{I}, \text{ for any  subset }I\subset\overline{1,n}.$
\end{theorem}

{\it Proof.}
We may assume that $t=1$ since this simplification does not hide any essential part of the argument. Indeed, note that for any $i\in\overline{1,n}$ there is a natural identification $M^t=M_{\widehat i}\bar\otimes M_i^ t$. 
Assume by contradiction that there exists $i$ such that $f(i)\neq i.$
Using the assumption we can apply Theorem \ref{Cor:ultrapower} and derive that for any $j\in\overline{1,n}$ there exists a  subgroup $\Sigma_j<\Gamma$ such that 
\begin{equation}\label{t3}
    M_{\widehat j}\prec_M^s L(\Sigma_j) \text{ and } L(\Sigma_j) \prec_M M_{\widehat {f(j)}}. 
\end{equation}
Here, we also used \cite[Lemma 2.4(2)]{DHI16}.
Note that \eqref{t3} shows that $L(\Sigma_i)\prec_M M_{\widehat{f(i)}}$ and $M_{\widehat{f(i)}}\prec_M^s L(\Sigma_{f(i)})$. Using \cite[Lemma 3.7]{Va08} we deduce that 
$L(\Sigma_i)\prec_M L(\Sigma_{f(i)})$. We can apply \cite[Lemma 2.2]{CI17} and derive that there exist $g\in\Gamma$ and a finite index subgroup $\Sigma_i^0<\Sigma_i$ such that $g\Sigma_i^0g^{-1}\subset \Sigma_{f(i)}$. Using \eqref{t3} for $j=f(i)$ we get that $L(\Sigma_i^0)\prec_M M_{\widehat {f^2(i)}}$. Since $[\Sigma_i:\Sigma_i^0]<\infty$, it follows that 
\begin{equation}\label{t4}
L(\Sigma_i)\prec_M M_{\widehat {f^2(i)}}.
\end{equation}
Following an idea from \cite[Section 3]{CD-AD20}, we consider a Bernoulli action $\Gamma\car D$ with abelian base and let $\mathcal M=D\rtimes\Gamma$. Let $\Psi:\mathcal M\to\mathcal M\bar\otimes M$ be the $*$-homomorphism given by $\Psi(du_g)=du_g\otimes u_g$, for all $d\in D$ and $g\in\Gamma$. From \eqref{t3} we have 
that $M_{\widehat i}\prec_M^s L(\Sigma_i)$, which gives $M_{\widehat i}\prec_{\mathcal M}^{s,s'} D\rtimes\Sigma_i$. By applying \cite[Lemma 2.3]{Dr19}, we further derive that  $\Psi(M_{\widehat i})\prec_{\mathcal M\bar\otimes M}^{s,s'} \mathcal M\bar\otimes L(\Sigma_i)$. Next, by using \cite[Lemma 2.4]{Dr19} and the relation \eqref{t4} we obtain that
\begin{equation}\label{t5}
\Psi(M_{\widehat i})\prec_{\mathcal M\bar\otimes M}^{s} \mathcal M\bar\otimes M_{\widehat {f^2(i)}}.
\end{equation}
Note that $i\in \widehat{ f(i)}.$ Using our assumption, we obtain that  $\Delta(M_{i})\prec^s_{M\bar\otimes M}M\bar\otimes M_{\widehat {f^2(i)}}$. Hence, by \cite[Remark 2.2]{DHI16} we get that
\begin{equation}\label{t6}
    \Psi(M_{i})\prec^s_{\mathcal M\bar\otimes M}M\bar\otimes M_{\widehat {f^2(i)}}.
\end{equation}
Next, by using \cite[Lemma 2.6]{Is19}, \eqref{t5} and \eqref{t6}, we derive that $\Psi(M)\prec^s_{\mathcal M\bar\otimes M}\mathcal M\bar\otimes  M_{\widehat {f^2(i)}}$. Since $\Psi(x)=x\otimes 1$, for any $x\in D$, we deduce that  $\Psi(\mathcal M)\prec^s_{\mathcal M\bar\otimes M}\mathcal M\bar\otimes  M_{\widehat {f^2(i)}}$. By \cite[Lemma 10.2]{IPV10} we obtain that $M\prec_{\mathcal M} M_{\widehat {f^2(i)}}$. Since $\mathcal {QN}_{\mathcal M}^{(1)}(M)=M$, we derive that $M\prec_{ M} M_{\widehat {f^2(i)}}$, which implies that $M_{f^2(i)}$ is not diffuse, contradiction. 

Hence, $\Delta(M_{\widehat j})\prec^s_{M\bar\otimes M}M\bar\otimes M_{\widehat j}, \text{ for any } j\in\overline{1,n}.$ By applying \cite[Lemma 2.6(2)]{DHI16}, the conclusion of the theorem follows.
\hfill$\blacksquare$

\section{Proofs of the main results}

\subsection{Proof of Theorem \ref{A}}

Before proceeding to the proof of Theorem \ref{A}, we need the following result.

\begin{theorem}\label{P:onebyone}
Let $M=M_1\bar\otimes\dots \bar\otimes M_n$ be the tensor product of $n\ge 1$ II$_1$ factors from $\mathscr M\cup\mathscr M_{wr}$.
Assume $\Gamma$ is a countable icc group such that $M^t=L(\Gamma)$ for some $t>0$.
Denote by $\Delta:M^t\to M^t\bar\otimes M^t$ the $*$-homomorphism given by $\Delta(u_g)=u_g\otimes u_g$, for any $g\in\Gamma$.

Then for any $i\in\overline{1,n}$, there exists $f(i)\in\overline{1,n}$ such that 
$\Delta(M_{\widehat i}^t)\prec^s_{M^t\bar\otimes M^t} M^t\bar\otimes M^t_{\widehat {f(i)}}.$
\end{theorem}

{\it Proof.} Let $(\tilde M_i, (\alpha^i_t)_{t\in\mathbb R})$ be an s-malleable deformation of $M_i$ given by the fact that $M_i$ belongs to $\mathscr M\cup\mathscr M_{wr}$, see Remark \ref{remarkbernoullimalleable}(1). If $M_i\in\mathscr{M}$, we let $A_i$ be given as in the definition of Class $\mathscr M$. 
We naturally extend $\alpha^i_t$ to an automorphism $\alpha_t^i\in\text{Aut}(M_{\widehat i}\bar\otimes \tilde M_i)$. We may assume without loss of generality that $t=1$ since this simplification does not hide any essential part of the argument; note that for any $i\in\overline{1,n}$ there is a natural identification $M^t=M_{\widehat i}\bar\otimes M_i^ t$. 

We first prove that there exists a function $f\in\overline{1,n}\to\overline{1,n}$ such that
\begin{equation}\label{t1}
\Delta(M_{j}) \text{ is non-amenable relative to } M\bar\otimes M_{\widehat {f(j)}} \text{ for any }j\in \overline{1,n}.    
\end{equation}
If \eqref{t1} does not hold, then there exists $j\in\overline{1,n}$ such that $\Delta(M_j) \text{ is amenable relative to } M\bar\otimes M_{\hat k},$ for any $k\in\overline{1,n}$. By \cite[Proposition 2.7]{PV11}, it follows that $\Delta(M_j) \text{ is amenable relative to } M\bar\otimes 1.$ This further implies by \cite[Proposition 7.2(4)]{IPV10} that $M_j$ is amenable, contradiction. Therefore, there exists  a function $f$ that satisfies \eqref{t1}.
Next, we show that
\begin{equation}\label{ttt}
    \Delta(M_{\widehat j})\prec_{M\bar\otimes M}M\bar\otimes M_{\widehat {f(j)}}, \text{ for any } j\in\overline{1,n}. 
\end{equation}

To show this, fix $j\in\overline{1,n}$ and notice first that $\mathcal{N}_{M\bar\otimes M}(\Delta(M_j))'\cap M\bar\otimes M=\mathbb C 1$ since $\Gamma$ is icc. Thus, using \cite[Lemma 2.6(2)]{DHI16}, we obtain from \eqref{t1} that
\begin{equation}\label{t14}
\Delta(M_j) \text{ is strongly non-amenable relative to }M\bar\otimes M_{\widehat {f(j)}}.    
\end{equation}

{\bf Case 1.} $M_{f(j)}\in\mathscr M.$

In this case it follows directly from \eqref{t14} that
$\text{id}\otimes\alpha_t^{f(j)}\to \text{id}$ uniformly on $(\Delta(M_{\widehat j}))_1$.
Assume by contradiction that $\Delta(M_{\widehat j})\nprec_{M\bar\otimes M}M\bar\otimes M_{\widehat {f(j)}}$. If $\Delta(M_{\widehat j})\nprec_{M\bar\otimes M}M\bar\otimes (M_{\widehat {f(j)}}\bar\otimes A_{f(j)})$, we obtain from Theorem \ref{Th:deform}(1) that $\text{id}\otimes \alpha_t^{f(j)}\to \text{id}$ uniformly on $(\Delta(M))_1$. This shows that $\alpha_t^{f(j)}\to \text{id}$ uniformly on $(M)_1$, and hence on $M_{f(j)}$, contradiction. Therefore, $\Delta(M_{\widehat j})\prec_{M\bar\otimes M}M\bar\otimes (M_{\widehat {f(j)}}\bar\otimes A_{f(j)})$. Using condition (2.ii) from Class $\mathscr M$, it follows that \eqref{ttt} holds in  this case.

{\bf Case 2.} $M_{f(j)}\in\mathscr M_{wr}$.

In this case we can write $M_{f(j)}=B_{f(j)}\rtimes\Lambda_{f(j)}$ as the von Neumann algebra of a generalized Bernoulli action as in Class $\mathscr M_{wr}$. By \eqref{t14} and \cite[Theorem 3.1]{BV12} we have that
$\text{id}\otimes\alpha_t^{f(j)}\to \text{id}$ uniformly on $(\Delta(M_{\widehat j}))_1$.
Using \cite[Theorem 4.2]{IPV10} we obtain that 
$\Delta(M_{\widehat j})\prec_{M\bar\otimes M}M\bar\otimes M_{\widehat {f(j)}}$ 
or $\Delta(M)\prec_{M\bar\otimes M} M\bar\otimes (M_{\widehat {f(j)}}\bar\otimes L(\Lambda_{f(j)}))$ or $\Delta(M) \prec_{M\bar\otimes M} M\bar\otimes (M_{\widehat {f(j)}}\bar\otimes (B_{f(j)}\rtimes\theta_{f(j)}))$, where $\theta_{f(j)}<\Lambda_{f(j)}$ is an infinite index subgroup. From \cite[Proposition 7.2]{IPV10} the last two possibilities give a contradiction. This shows that \eqref{ttt} holds.

Finally, by using \cite[Lemma 2.4(2)]{DHI16}, we end the proof of the theorem.
 \hfill$\blacksquare$

We are now ready to prove the following result which is a generalization of
Theorem \ref{A}.

\begin{theorem}\label{Th:general}
Let $M=M_1\bar\otimes\dots \bar\otimes M_n$ be the tensor product of $n\ge 1$ II$_1$ factors from $\mathscr M\cup\mathscr M_{wr}$.
Assume $\Gamma$ is a countable icc group such that $M^t=L(\Gamma)$ for some $t>0$.

Then there exist a product decomposition $\Gamma=\Gamma_1\times\dots\times \Gamma_n$, a unitary $u\in M$ and some positive numbers $t_1,\dots,t_n$ with $t_1\cdot\dots\cdot t_n=t$ such that $M^{t_i}_i=uL(\Gamma_i)u^*$, for any $i\in \overline{1,n}$.
\end{theorem}

{\it Proof.} 
Let $(\tilde M_i, (\alpha^i_t)_{t\in\mathbb R})$ be an s-malleable deformation of $M_i$ given by the fact that $M_i$ belongs to Class $\mathscr M\cup\mathscr M_{wr}$, see Remark \ref{remarkbernoullimalleable}(1). If $M_i\in\mathscr{M}$, we let $A_i$ be given as in the definition of Class $\mathscr M$. 
We naturally extend $\alpha^i_t$ to an automorphism $\alpha_t^i\in\text{Aut}(M_{\widehat i}\bar\otimes \tilde M_i)$. Next, we may assume that $t=1$ since this simplification does not hide any essential part of the argument.

Following \cite{PV09}, we denote by $\Delta:M\to M\bar\otimes M$ the $*$-homomorphism given by $\Delta(u_g)=u_g\otimes u_g$, for any $g\in\Gamma$. By Theorem \ref{P:onebyone} and Theorem \ref{Th:generalonebyone}, we have that  
\begin{equation}\label{q1}
\Delta(M_{I})\prec^s_{M\bar\otimes M} M\bar\otimes M_{I}, \text{ for any  subset }I\subset\overline{1,n}.    
\end{equation}

Our goal is to prove the following claim.

{\bf Claim.} There exist a subgroup $\Sigma<\Gamma$ and a non-zero projection $f\in L(\Sigma)'\cap M$ such that 
$
M_{\widehat n}\prec_M L(\Sigma)f \text{ and } L(\Sigma)f\prec_M^s M_{\widehat n}. 
$

{\it Proof of Claim.}
Using Theorem \ref{Cor:ultrapower}, we can find a subgroup $\Sigma<\Gamma$ with non-amenable centralizer $C_\Gamma(\Sigma)$ such that $M_{\hat n}\prec_M L(\Sigma) \text{ and } L(\Sigma)\prec_M M_{\hat n}.$ By using \cite[Lemma 2.4(4)]{DHI16} there exists
a non-zero projection $f\in L(\Sigma C_\Gamma(\Sigma))'\cap M$ such that
\begin{equation}\label{q2}
M_{\hat n}\prec_M ^{s'}L(\Sigma)f \text{ and } L(\Sigma)\prec_M M_{\hat n}.    
\end{equation}

Next, by using \cite[Lemma 3.5]{Va08} and \cite[Lemma 2.6(3)]{DHI16} we deduce from \eqref{q2} that $L(C_\Gamma(\Sigma))f$ is amenable relative to $M_{n}$. Since $C_\Gamma(\Sigma)$ is non-amenable, we derive from \cite[Proposition 2.7]{PV11} that $L(C_\Gamma(\Sigma))f$ is non-amenable relative to $M_{\widehat n}$. Using \cite[Lemma 2.4]{DHI16} there exists a non-zero projection $f_1\in L(\Sigma C_\Gamma(\Sigma))'\cap M$ with $f_1\leq f$ such that 
\begin{equation}\label{q3}
    L(C_\Gamma(\Sigma))f_1 \text{ is strongly non-amenable relative to } M_{\widehat n}.
\end{equation}

Next, we note that 
\begin{equation}\label{extra}
\alpha_t^n\to \text{id} \text{ uniformly on }(L(\Sigma)f_1)_1. 
\end{equation}
Indeed, if $M_n$ belongs to $\mathscr{M}$, this follows immediately, while if $M$ belongs to $\mathscr M_{wr}$, this follows from \cite[Theorem 3.1]{BV12}.
Next, remark that the claim would follow if $L(\Sigma)f_1\prec_M M_{\widehat n}$ because we could use \cite[Lemma 2.4]{DHI16} and derive that there exists a non-zero projection $f_0\in L(\Sigma)'\cap M$ with $f_0\leq f_1$ such that $L(\Sigma)f_0\prec^s_M M_{\hat n}$.

Hence, we assume by contradiction that $L(\Sigma)f_1\nprec_M M_{\widehat n}$. 

Following \cite[Section 4]{CdSS15}, we let  $\Omega=\{g\in\Gamma|\mathcal O_{\Sigma}(g) \text{ is finite}\}$, where $\mathcal O_\Sigma(g)=\{hgh^{-1}|\; h\in\Sigma\}$, and notice that $L(\Sigma)'\cap M\subset L(\Omega)$. We continue by proving the following.

{\bf Subclaim.} 
There exists a non-zero projection $f_2\in L(\Omega)'\cap M$ such that 
\[\alpha_t^n\to \text{id}  \text{ uniformly on }(L(\Omega)f_2)_1.\]

{\it Proof of subclaim.} We split the proof of the subclaim in two parts.

{\bf Case 1.} $M_n$ belongs to Class $\mathscr M$.

Note that there exists a sequence of increasing subgroups $\Omega_i<\Omega$ such that $\Omega=\vee_{i\ge 1}\Omega_i$ and by letting $\Sigma_i=C_\Sigma (\Omega_i)$, we have $[\Sigma:\Sigma_i]<\infty$, for any $i\ge 1$. Indeed, let $\{\mathcal O_i\}_{i\ge 1}$ be a countable enumeration of all the finite orbits of the action by conjugation of $\Sigma$ on $\Gamma$. Notice that $\Omega_i:=\langle \cup_{j=1}^i \mathcal O_j \rangle<\Omega$, $i\ge 1$, is an ascending sequence of subgroups with $\Omega=\vee_{i\ge 1}\Omega_i$. Since $\cup_{j=1}^i \mathcal O_j \subset\Omega$ is a finite set, it follows that $\Sigma_i:=\cap_{g\in \cup_{j=1}^i \mathcal O_j } C_\Sigma(g)=C_\Sigma(\Omega_i)$ is a decreasing sequence of finite index subgroups of $\Sigma$.

Next, note that $f_1\in L(\Sigma_i)'\cap M$ and $L(\Sigma_i)f_1\nprec_M M_{\hat n}$ since $[\Sigma:\Sigma_i]<\infty$, for any $i\ge 1$. Moreover, we note that $L(\Sigma_i)f_1\nprec_M M_{\hat n}\bar\otimes A_n$. Otherwise, by using condition (2) in Class $\mathscr M$, we get that there exists a non-zero projection $f_2\in (L(\Sigma_i)'\cap M)'\cap M\subset L(C_{\Gamma}(\Sigma))'\cap M$ with $f_2\leq f_1$ such that $(L(\Sigma_i)'\cap M)f_2$ is amenable relative to $M_{\widehat n}$, which contradicts \eqref{q3}. 

Therefore, since $M\subset M_{\widehat n}\bar\otimes \tilde M_n$ is mixing relative to $M_{\widehat n}\bar\otimes A_n$ we obtain from \eqref{extra} and Theorem \ref{Th:deform}(1) that  
\begin{equation}\label{q4}
\alpha_t^n\to \text{id}  \text{ uniformly on } (f_1(L(\Sigma_i)\vee (L(\Sigma_i)'\cap M))f_1)_1, \text{ for any }i\ge 1.    
\end{equation}
Next, since $A_n$ is amenable, we notice that \eqref{q3} together with \cite[Theorem 2.6(3)]{DHI16} imply that $f_1(L(\Sigma_i)'\cap M))f_1\nprec_M M_{\widehat n}\bar\otimes A_n$ for any $i\ge 1$. We can therefore combine \eqref{q4} with Theorem \ref{Th:deform}(2) and deduce that
\[
\alpha_t^n\to \text{id}  \text{ uniformly on } (f_1(\bigvee_{i\ge 1}(L(\Sigma_i)'\cap M))f_1)_1.
\]
We can apply Proposition \ref{P:central} and derive that if we denote by $f_2$ the central support of $f_1$ in $\bigvee_{i\ge 1}(L(\Sigma_i)'\cap M)$, we further obtain from that
\[
\alpha_t^n\to \text{id}  \text{ uniformly on } ((\bigvee_{i\ge 1}(L(\Sigma_i)'\cap M))f_2)_1.
\]
Since $L(\Omega)\subset \bigvee_{i\ge 1}L(\Sigma_i)'\cap M$ and $f_2\in \mathcal Z(\bigvee_{i\ge 1}(L(\Sigma_i)'\cap M))$, it follows that the subclaim is proven in this case.

{\bf Case 2.} $M_n$ belongs to Class $\mathscr M_{wr}$.

We can write $M_n=B_n\rtimes\Lambda_n$ where $\Lambda_n\car B_n$ is a generalized Bernoulli action with amenable base as given by Class $\mathscr M_{wr}$. For proving the subclaim, we follow a slightly different approach than the one used in Case 1. Using Popa's compression formulas for quasi-normalizers \cite{Po03}, we have 
\begin{equation}\label{qn01}
\mathcal {QN}_{f_1Mf_1}(L(\Sigma)f_1)= f_1 \mathcal {QN}_{M}(L(\Sigma))f_1.    
\end{equation}
By applying \cite[Theorem 4.2]{IPV10} to \eqref{extra}, we derive that (i) $L(\Sigma)f_1\prec_M M_{\widehat n}$ or (ii) $L(QN_{\Gamma}(\Sigma))f_1$ $\prec_M M_{\widehat n}\bar\otimes (B_n\rtimes \theta_n)$, where $\theta_n<\Lambda_n$ is an infinite index subgroup, or (iii) there exists a partial isometry $w\in M$ with $ww^*=f_1$ and $w^* \mathcal {QN}_{f_1Mf_1}(L(\Sigma)f_1)'' w \subset M_{\widehat n}\bar\otimes L(\Lambda_n)$. 
Option (i) is not possible since we assumed by contradiction that $L(\Sigma)f_1\nprec_M M_{\widehat n}$. We now show that option (ii) leads to a contradiction as well. Note that by passing to relative commutants in \eqref{q2}, we derive that $M_n\prec_M L(\Omega)$. Since $\Omega\subset QN_{\Gamma}(\Sigma)$, option (ii) implies that $L(\Omega)\prec_M M_{\widehat n}\bar\otimes (B_n\rtimes \theta_n)$. Combining all these with \eqref{q2}, we can apply Proposition \ref{augm} and derive that $[\Lambda_n:\theta_n]<\infty$, contradiction.


Next, note that option (iii) combined with \eqref{qn01} implies that 
\[
\alpha_t^n\to \text{id}  \text{ uniformly on } (f_1 \mathcal {QN}_{M}(L(\Sigma))''f_1)_1.
\]
Since $\Omega\subset QN_{\Gamma}(\Sigma)$,
the subclaim follows by using Proposition \ref{P:central}.
\hfill$\square$



Following the augmentation technique from \cite[Section 3]{CD-AD20}, we consider a Bernoulli action $\Gamma\car D$  with abelian base. Denote $\mathcal M=D\rtimes\Gamma$ and let $\Psi:\mathcal M\to \mathcal M\bar\otimes M$ be the $*$-homomorphism given by $\Psi(du_g)=du_g\otimes u_g$, for all $d\in D$ and $g\in\Gamma$.
Next, by passing to relative commutants in \eqref{q2} we obtain that $M_n\prec_M L(\Omega)$, which implies that $M_n\prec_{\mathcal M}^{s'} D\rtimes\Omega$. This shows using \cite[Lemma 2.3]{Dr19} that $\Psi(M_n)\prec_{\mathcal M\bar\otimes M}^{s'}\mathcal M\bar\otimes L(\Omega)$. In particular, $\Psi(M_n)\prec_{\mathcal M\bar\otimes M}\mathcal M\bar\otimes L(\Omega)f_2$.

Hence, there exist some projections $p\in M_n,q\in \mathcal M\bar\otimes L(\Omega)f_2$, a non-zero partial isometry $w\in q(\mathcal M\bar\otimes M)\Psi(p)$ and a $*$-homomorphism $\theta: \Psi(pM_np)\to q(\mathcal M\bar\otimes L(\Omega)f_2)q$ such that $\theta(x)w=wx$, for any $x\in \Psi (pM_np).$ Let $\tilde p=\Psi(p)$ and note that $w^*w\in \Psi(pM_np)'\cap \tilde p(\mathcal M\bar\otimes M)\tilde p$. Therefore, by using the multiplicativity of $\text{id}\otimes\alpha_t^n$ and its pointwise $\|\cdot\|_2$-convergence to the identity, the subclaim implies that
\begin{equation}\label{q6}
\text{id}\otimes\alpha_t^n\to \text{id} \text{ uniformly on } (\Psi(pM_np)w^*w)_1.    
\end{equation}

The claim will be proven by considering again two separate cases.

{\bf Case 1.} $M_n$ belongs to Class $\mathscr M$.

If $\Psi(pM_np)w^*w\nprec_{\mathcal M\bar\otimes M}\mathcal M\bar\otimes (M_{\hat n}\bar\otimes A_n)$, then Theorem \ref{Th:deform}(1) and \eqref{q6} give
\[
\text{id}\otimes\alpha_t^n\to \text{id} \text{ uniformly on } (w^*w(\Psi(pM_np)\vee \Psi(pM_np)'\cap \tilde p(\mathcal M\bar\otimes M)\tilde p)w^*w)_1.
\]
Let $\tilde p_1$ be the central support of $w^*w$ in $\Psi(M_n)\vee \Psi(M_n)'\cap  (\mathcal M\bar\otimes M)$ and note that $\tilde p_1\in \Psi(M)'\cap (\mathcal M\bar\otimes M)=\mathbb C 1$ since $\Gamma$ is icc. By using Proposition \eqref{P:central}, we further obtain that $\text{id}\otimes\alpha_t^n\to \text{id} \text{ uniformly on  } (\Psi(M))_1,$ which shows that $\alpha_t^n\to \text{id}  \text{ uniformly on } (M)_1,$ contradiction.

Hence, $\Psi(pM_np)w^*w\prec_{\mathcal M\bar\otimes M}\mathcal M\bar\otimes (M_{\hat n}\bar\otimes A_n)$, and therefore $\Psi(M_n)\prec_{\mathcal M\bar\otimes M}^s\mathcal M\bar\otimes (M_{\hat n}\bar\otimes A_n)$ by \cite[Lemma 2.4(2)]{DHI16}. Relation \eqref{q1} implies in particular that $\Psi(M_n)\prec^s_{\mathcal M\bar\otimes M}\mathcal M\bar\otimes M_{n}$. By using \cite[Lemma 2.8(2)]{DHI16}, we derive that $\Psi(M_n)\prec_{\mathcal M\bar\otimes M}^s\mathcal M\bar\otimes A_n$, which implies by \cite[Lemma 10.2]{IPV10} that $M_n$ is amenable, contradiction. Hence, the claim is proven in this case.

{\bf Case 2.} $M_n$ belongs to Class $\mathscr M_{wr}$.

As before, we assume $M_n=B_n\rtimes\Lambda_n$ where $\Lambda_n\car C_n^{I_n}=:B_n$ is a generalized Bernoulli action with amenable base as given by Class $\mathscr M_{wr}$. If $\Psi(pM_np)w^*w\nprec_{\mathcal M\bar\otimes M}\mathcal M\bar\otimes M_{\hat n}$,
we derive from  \cite[Theorem 4.2]{IPV10} and \eqref{q6} that $\Psi(M)\prec_{\mathcal M\bar\otimes M} \mathcal M\bar\otimes (B_n\rtimes\theta_n)$, where $\theta_n<\Lambda_n$ is an infinite index subgroup, or $\Psi(M)\prec_{\mathcal M\bar\otimes M} \mathcal M\bar\otimes L(\Lambda_n)$, and therefore, $\Psi(\mathcal M)\prec_{\mathcal M\bar\otimes M} \mathcal M\bar\otimes (B_n\rtimes\theta_n)$ or $\Psi(\mathcal M)\prec_{\mathcal M\bar\otimes M} \mathcal M\bar\otimes L(\Lambda_n)$, respectively. By using \cite[Lemma 10.2]{IPV10} it is easy to see that we would get a contradiction. Hence, $\Psi(pM_np)w^*w\prec_{\mathcal M\bar\otimes M}\mathcal M\bar\otimes M_{\hat n}$, and by proceeding as in the last paragraph of Case 1 above, we obtain that 
$\Psi(M_n)\prec_{\mathcal M\bar\otimes M}^s\mathcal M\otimes 1$, which is again a contradiction.
\hfill$\square$

Finally, note that \eqref{q1} implies in particular that $\Delta(M_1)\prec_{M\bar\otimes M}M\bar\otimes M_1$. Using Theorem \ref{Cor:ultrapower}, we obtain a subgroup $\Theta<\Gamma$ such that $M_1\prec_M L(\Theta)$ and $L(\Theta)\prec_M M_1$. In combination with the claim, it follows from Theorem \ref{Th:splits} that there exist a product decomposition $\Gamma=\Gamma_1^{n-1}\times\Gamma_n$, a decomposition $M=M_{\widehat n}^s\bar\otimes M_n^{1/s}$ for some $s> 0$ and a unitary $u\in\mathcal U(M)$ such that $uM_{\widehat n}^su^*=L(\Gamma_1^{n-1})$ and $uM_n^{1/s}u^*=L(\Gamma_n)$. 
Therefore, we obtain the conclusion of the theorem by a simple induction argument.
\hfill$\blacksquare$

\subsection{Proof of Corollary \ref{B}}

Let $\theta:L(\Gamma)^t\to L(\Lambda)$ be a $*$-ismorphism where $\Lambda$ is any countable group and $t>0$. By Theorem \ref{A}, there exist a product decomposition $\Lambda=\Lambda_1\times\dots\times\Lambda_n$, some positive numbers $t_1,\dots,t_n>0$ with $t_1\cdots t_n=t$ and a unitary $w\in L(\Lambda)$ such that $\theta(L(\Gamma_i)^{t_i})=w L(\Lambda_i)w^*$, for any $i\in\overline{1,n}.$ Since $\Gamma_i$ is W$^*$-superrigid, it follows that $t_i=1$ and there exist a group isomorphism $\delta_i:\Gamma_i\to\Lambda_i$, a unitary $w_i\in L(\Lambda_i)$ and a character $\omega_i:\Gamma_i\to \mathbb T$ such that $\theta(u_g)=\omega_i(g) w_i v_{\delta(g)}w_i^*$, for all $i\in\overline{1,n}$ and $g\in\Gamma_i$. Hence, $t=1$ and by letting  $\omega = \prod_{i=1}^{n}{\omega_{i}}$,  $\delta = \prod_{i=1}^{n}{\delta_{i}}$ and $u=\prod^n_{i=1}w_i$, we get the desired conclusion. Finally, we notice that any group from $\mathcal {IPV}$ is W$^*$-superrigid by \cite[Theorem 8.3]{IPV10}.
\hfill$\blacksquare$

\subsection{Proof of Corollary \ref{C}}
On one hand, note that if $\Gamma_0$ is an icc non-amenable countable group with $\beta^{(2)}_1(\Gamma_0)>0$, then its amenable radical is trivial, i.e. any normal amenable subgroup of $\Gamma_0$ is trivial. On the other hand, note that any group from $\mathcal A_2$ has trivial amenable radical by using \cite[Proposition 6.3]{CD-AD20}.

Therefore, $\Gamma$ has trivial amenable radical, so by using \cite[Theorem 1.3]{BKKO14}, it follows that $C^*_r(\Gamma)$ has a unique trace. This implies that any $*$-isomorphism $\theta:C^*_r(\Gamma)\to C^*_r(\Lambda)$ extends to a $*$-isomorphism $\theta:L(\Gamma)\to L(\Lambda)$. The conclusion now follows from Theorem \ref{A}. 
\hfill$\blacksquare$

\end{document}